\documentstyle[epsf,amssymb,amsfonts,12pt]{amsart}

\newcommand{\calf}{{\cal F}}
\newcommand{\Chow}{\mbox{\it Chow}\,}

\newcommand{\DOT}{\setlength{\unitlength}{1pt}\begin{picture}(2.5,2)(1,1)
\put(1,2){\circle*{2}}\end{picture}}

\newcommand{\Fdot}{{F\!_{\DOT}}}

\newcommand{\Fl}{\mbox{\bf Fl}}

\newcommand{\G}{\mbox{\bf G}}
\newcommand{\Gn}{{\mbox{\bf G}_1\Pn}}
\newcommand{\GI}{{\mbox{\bf G}_1}}

\newcommand{\Hom}{\mbox{Hom}}
\newcommand{\hfl}[1]{\left\lfloor \frac{#1}{2} \right\rfloor}

\newcommand{\Ndot}{{N_{\DOT}}}

\newcommand{\Pp}{\mbox{\bf P}}
\newcommand{\Pn}{{\mbox{\bf  P}^n}}
\newcommand{\PI}{{\mbox{\bf P}^1}}

\newcommand{\QED}{
\setlength{\unitlength}{1.0pt}%
\begin{picture}(7.5,7.5)
\put(0,-2.5){\rule{2.5pt}{5pt}}
\put(0,2.5){\rule{5pt}{2.5pt}}
\put(0,5){\rule{7.5pt}{2.5pt}}
\end{picture}\vspace{10pt}}

\newcommand{\R}{\mbox{\bf R}}
\newcommand{\Rsc}{\mbox{\scriptsize\bf R}}

\newcommand{\Span}[1]{\langle #1 \rangle}

\newsavebox{\Shading}

\sbox{\Shading}{\setlength{\unitlength}{1.3pt}%
\begin{picture}(30,20)

\put( 2,11.8){\bf 1}
\put(12,11.8){\bf 2}
\put(22,11.8){\bf 3}
\put( 2, 1.8){\bf 4}
\end{picture}}

\begin{document}

\title[Real Enumerative geometry]{Explicit Enumerative Geometry for the
Real Grassmannian of Lines in Projective Space}

\author{Frank Sottile}

\address{
        Department of Mathematics\\
        University of Toronto\\
        100 St. George Street\\
	Toronto, Ontario  M5S 1A1\\
	Canada\\
	(416) 978-4031}
\email{sottile@@math.toronto.edu}
\date{31 October 1995}
\thanks{Research supported in part by NSERC grant \# OGP0170279}

\begin{abstract}
For any collection of Schubert conditions on lines in
projective space which generically determine a finite number of lines,
we show there exist real generic conditions determining the expected
number of real lines.
This extends the classical Schubert calculus of enumerative geometry
for the Grassmann variety of lines in projective space from the complex
realm to the real.
Our main tool is an explicit description of
rational equivalences which also constitutes a novel determination of the
Chow rings of these Grassmann varieties.
The combinatorics of these rational equivalences
suggests a non-commutative associative
product on the free abelian group on Young tableaux.
We conclude by considering some enumerative problems over other
fields.
\end{abstract}

\maketitle

\section{Introduction}
A basic problem in algebraic geometry is to describe the common
zeroes of a set of polynomials.
This is more difficult over  non-algebraically closed fields.
For systems of polynomials with few monomials on a complex torus,
Khovanskii~\cite{Khovanskii_fewnomials} showed
that the real zeroes are at most a small fraction of
the complex zeroes.
Fulton (\cite{Fulton_introduction_intersection}, \S7.2) asked
how many solutions to a problem of
enumerative geometry can be real; in particular, how many of the
3264 conics tangent to five general real conics can be real.
He later showed that all, in fact, can be real.
Recently, this was independently rediscovered by Ronga, Tognoli, and
Vust~\cite{Ronga_Tognoli_Vust}.
Robert Spesier suggested the classical Schubert calculus of enumerative
geometry would be a
good testing ground for these questions.
For any problem of enumerating lines in ${\bf P}^n$ incident on
real linear subspaces in general position, we show
that all solutions can be real.
\medskip

A flag and a partition $\lambda = (\alpha,\beta)$ determine a Schubert
subvariety of the Grassmannian of lines in ${\bf P}^n$
of type $\lambda$, which has codimension
$|\lambda| = \alpha+\beta$.
Any generically transverse intersection of
Schubert varieties  is rationally equivalent
to a sum of Schubert varieties.
The classical Schubert calculus gives algorithms for determining
how many of each type.

For  partitions $\lambda^1,\ldots,\lambda^m$,
we describe a cycle $\Omega({\cal T})$ (depending upon $\lambda^1,
\ldots,\lambda^m$) which is
a sum of distinct Schubert varieties.
Let ${\cal G}= {\cal G}(\lambda^1,\ldots,\lambda^m)$
be the set of points of the Chow variety
representing  generically transverse intersections
of Schubert varieties of types $\lambda^1,\ldots,\lambda^m$.
In \S 4, we show that ${\cal G}$ is unirational.
Cycles represented by the points of a rational curve on the Chow
variety are rationally equivalent.
In \S 3, we show
\medskip

\noindent{\bf Theorem A.} \ {\em
Let $X\in {\cal G}$.
Then there is a chain of rational curves between
$X$ and the cycle $\Omega({\cal T})$.
Furthermore, these curves may be explicitly
described and each lies in the
Zariski closure of ${\cal G}$.
In particular, the point representing
$\Omega({\cal T})$ is in the
Zariski closure of $\,{\cal G}$.}
\medskip

The proof of Theorem A constitutes an explicitly geometric
determination of the Schubert
calculus of enumerative geometry for lines in ${\bf P}^n$.
In fact, it shows these `Schubert-type' enumerative problems
may be solved {\em without} reference to the Chow ring, a traditional tool
in enumerative geometry.
We use it to compute products in the Chow ring.
Let $\sigma_\lambda$ be the rational equivalence class
of a Schubert variety of type $\lambda$.
\medskip

\noindent{\bf Theorem B. }
{\em
Let $c^\lambda$ be the number  of components of $
\,\Omega({\cal T})$ of type
$\lambda$.
Then
$$
\prod_{i=1}^m \sigma_{\lambda^i} = \sum_\lambda c^\lambda \sigma_\lambda.
$$
}

Thus we derive the structure of these
Chow rings in a strong sense:
All products among classes from the Schubert basis
are expressed as linear combinations of basis elements
and these expressions are obtained by exhibiting
rational equivalences between a generically transverse intersection
of Schubert varieties and the cycle $\Omega({\cal T})$.
We believe this is the first non-trivial explicit description of
rational equivalences giving all products among a set
generators of the Chow group for any variety.

When $k= \R$, we show
\medskip

\noindent{\bf Theorem C.} \ {\em
Let $\lambda^1,\ldots,\lambda^m$ be partitions
with $|\lambda^1|+\cdots+|\lambda^m|$ equal to the dimension
of the Grassmannian of lines in $\mbox{\bf P}^n$.
Then there exists a nonempty classically open
subset  in the product of $m$ real flag manifolds whose
corresponding Schubert varieties meet transversally,
with all points of intersection real.
}
\medskip

To the best of our knowledge, this is the first result
showing that a large class of non-trivial enumerative problems can have all
of their solutions real.

The construction of the cycles  $\Omega({\cal T})$  and
rational curves  of Theorem~A use
a `calculus of tableaux' outlined
in \S\ref{sec:calculus_of_tableaux} and extended in \S6,
where we define a non-commutative associative algebra with additive basis
the set of Young tableaux.
This algebra has surjections to the Chow rings of Grassmann
varieties and the algebra of symmetric functions.
However, it differs fundamentally from the plactic algebra
of Lascoux and
Sch\"utzenberger~\cite{Lascoux_Schutzenberger_monoid_plactic}.

In \S7, we ask which enumerative
problems may be solved over which (finite) fields
and give the answer for two classes of Schubert-type
enumerative problems.
We also show how some of our constructions may be carried out over
finite fields.

The rational equivalences we construct are a modification of the
classical method of degeneration.
This method may fail when applied to more than
a few conditions; an intersection typically
becomes improper before the conditions become special enough to
completely determine the intersection.
Considering deformations of intersection cycles, rather
than of  conditions, an idea of
Chaivacci and
Escamilla-Castillo~\cite{Chiavacci_Escamilla-Castillo},
enables
us to deform generically transverse
intersections into sums of distinct Schubert varieties.

Theorem C is from our 1994
Ph.D. Thesis  from the University of Chicago,
written under the direction of William Fulton.
We would like to thank William Fulton for suggesting these problems,
for his thoughtful advice, and above all for introducing us
to algebraic geometry.

\section{Preliminaries}
Let $k$ be  an infinite field.
Varieties will be closed, reduced, projective
(not necessarily irreducible), and defined over $k$.
When $k={\R}$, let  $X(\R)$ be the points of $X$
with residue field $\bf R$.
We use the classical topology on $X(\R)$.

A subset $Y$ (not necessarily algebraic) of a variety
$X$ is {\em unirational} if $Y$ contains the image $U$ of a dense
open subset of affine space under an algebraic morphism and
 $Y\subset \overline{U}$.
A subset $Y\subset X(\R)$ is {\em real unirational} if
$Y$ contains the image $U$ of a dense open subset of $\R^n$ under a
real algebraic map and $Y\subset \overline{U}$.

Let $X$ be a smooth variety,
$U$ and $W$ subvarieties of $X$, and set $Z = U\cap W$.
Then $U$ and $W$ meet
{\em generically transversally}
if $U$ and $V$ meet transversally at the generic point of each component of
$Z$.
Then  $Z$ is generically reduced,
the fundamental cycle $[Z]$ of $Z$ is multiplicity free,
and in the Chow ring $A^*X$ of $X$:
$$
[U] \cdot [W] = [U\cap W] = [Z] = \sum_{i=1}^r\, [Z_i],
$$
where $Z_1,\ldots,Z_r$ are the irreducible components of $Z$.
\vspace{.5in}

\subsection{Chow Varieties}\label{sec:Chow}
Let $X$ be a projective variety.
The Chow variety $\Chow X$ is a projective variety parameterizing
positive cycles
on $X$.
Let $U$ be  a smooth variety and $W$ a subvariety
of $X\times U$ with equidimensional fibres over $U$.
Then there is an dense open subset $U'$ of $U$ such that
the association of a point $u$ of $U$ to the fundamental cycle of the
fibre $W_u$ determines a morphism $U'\rightarrow \Chow X$.
If $U$ is a smooth curve, then $U' = U$.
Moreover, if $X$, $U$, and $W$ are defined over $k$, then so
are $\Chow X$, $U'$, and the map $U' \rightarrow\Chow X$
(\cite{Samuel}, \S I.9).
Cycles represented all points on a rational curve
in  $\Chow X$ are rationally equivalent.
We will use  the same notation for a subscheme
of $X$, its fundamental cycle and the point representing that cycle
in $\Chow X$.
\vspace{.5in}

\subsection{Grassmannians and Schubert Subvarieties}
For $S\subset \Pn$, let $\Span{S}$ be its linear span.
For a vector space $V$, let ${\bf P} V$ be
the projective space of all one dimensional subspaces of $V$.
Suppose $K = {\Pp}U$ and $M={\Pp}W$.
Set
$\Hom(K,M) = \Hom(U,W)$, the space of linear maps from
$U$ to $W$.
If $K\subset M$, set $M/K = {\Pp}(W/U)$.
A complete flag $\Fdot$ is a collection of subspaces
$F_n\subset\cdots\subset F_1\subset F_0 = \Pn$,
where $\dim F_i = n-i$.
If $p>n$, set $F_p = \emptyset$.

For $0\leq s\leq n$, let $\G_s\Pn$ be the Grassmannian of $s$-dimensional
subspaces of $\Pn$.
We sometimes write
$\G^{n-s}\Pn$ for this variety.
Its dimension is $(n-s)(s+1)$.
For a partition $\lambda=(\alpha,\beta)$,
let $\Fl(\lambda)$
denote the variety of partial flags of {\em type} $\lambda$; those
$K\subset M$ with  $K$ a $(n-\alpha-1)$-plane
and $M$ a $(n-\beta)$-plane.

A partial flag $K\subset M$
determines a {\em Schubert variety} $\Omega(K,M)$;
those lines contained in $M$ which also meet $K$.
The {\em type} of $\Omega(K,M)$ is the type, $\lambda=(\alpha,\beta)$,
of $K\subset M$ and its codimension is $|\lambda| =\alpha+\beta$.
If $\alpha=\beta$, then $\Omega(K,M) = \G_1M$, the Grassmannian of lines
in $M$.
If $\beta=0$, so $M=\Pn$, then we write $\Omega_K$ for
this Schubert variety.

The tangent space to  $\ell \in \Gn$ is naturally
identified with the linear space $\Hom(\ell,\Pn/\ell)$.
It is not hard to verify the following Lemma, whose proof we omit.

\subsection{Lemma.}\label{lemma:one}
{\em
\begin{enumerate}
\item The smooth locus of $\Omega(K,M)$ consists of those
$\ell$ with $\ell \not\subset K$.
For such $\ell$,
$$
T_{\ell}\Omega(K,M) = \{ \phi\in \Hom(\ell,\Pn/\ell)\,|\,
\phi(\ell)\subset M/\ell\mbox{ and }
\phi(\ell\cap K)\in (K+\ell)/\ell\}.
$$
\item Let $K,M\subset {\bf P}^n$.
Then $\Omega_K\bigcap\G_1M =
\Omega(K\cap M,\, M)$,
and this is transverse at the
smooth points of $\,\Omega(K\cap M,\, M)$
if and only if $K$ and $M$ meet properly in $\Pn$.
\item Let   $K_i\subset M_i$, for $i=1,2$.
Then the intersection $\Omega(K_1,M_1)\bigcap
\Omega(K_2,M_2)$ is improper unless
$M_i$ meets both $K_j$ and $M_j$ properly for $i\neq j$.
\end{enumerate}
}
\medskip

An intersection of two Schubert varieties
may be  generically transverse and reducible:

\subsection{ Lemma.}\label{lemma:component_calculation}
{\em
Let $F,P,N,H$ be linear subspaces of ${\bf P}^n$ and suppose
$H$ is a hyperplane not
containing $P$ or $N$,
$F\subsetneq P\cap H$, and $N$ meets $F$ properly.
Set $L=N\cap H$.
Then $\Omega(F,P)$ and $\Omega_L$ meet generically transversally,
$$
\Omega(F,P)\bigcap\Omega_L \ =\
\Omega(N\cap F,P)    +\Omega(F,P\cap H)\bigcap\Omega_N,
$$
and the second component is itself a generically transverse intersection.
}
\medskip

\noindent{\bf Proof:}
The right hand side is a subset of the left, we show the other inclusion.
Let $\ell \in\Omega(F,P)\bigcap\Omega_L$.
If $\ell$ meets $L\cap F= N\cap F$, then
$\ell \in \Omega(N\cap F,P)$.
Otherwise, $\ell$ is spanned by its intersections with $F$ and $L$,
hence $\ell\subset \Span{F,L}\cap P \subset P\cap H$.
Since $N\cap P\cap H = L\cap P\cap H$, we see that
$\ell \in \Omega(F,P\cap H)\bigcap\Omega_N$.

Verifying these intersections are generically
transverse is left to the reader.
\QED

\subsection{A Calculus of Tableaux}\label{sec:calculus_of_tableaux}

The {\em Young diagram} of a partition $\lambda = (\alpha,\beta)$
is  a  two rowed array of boxes with
$\alpha$ boxes in the first row and $\beta$ in the second.
Note that $\alpha\geq \hfl{|\lambda|+1}\geq \hfl{|\lambda|}
\geq \beta$.
We make no distinction between a partition and its Young diagram.

A {\em Young tableau} $T$  of {\em shape} $\lambda$ is a
filling of the boxes of $\lambda$ with the integers $1,2,\ldots,|\lambda|$.
These integers increase left to right across each row and down each column.
Thus the $j$th entry in the second row of $T$ must be at
least $2j$.
Call $|\lambda|$ the {\em degree} of $T$, denoted $|T|$.
If $\alpha = \beta$, then $T$ is said to be {\em rectangular}.

Let $T$ be a tableau and $\alpha$ a positive integer.
Define $T* \alpha$ to be the set of all tableaux of degree
$|T|+\alpha$ whose first $|T|$ entries comprise $T$
and last $\alpha$ entries occur in distinct columns,
increasing from left to right.
For example:
$$
\setlength{\unitlength}{1.3pt}%
\begin{picture}(35,20)(0,8)
\thicklines
\put(0, 0){\line(1,0){10}}
\put(0,10){\line(1,0){30}}
\put(0,20){\line(1,0){30}}
\put( 0, 0){\line(0,1){20}}
\put(10, 0){\line(0,1){20}}
\put(20,10){\line(0,1){10}}
\put(30,10){\line(0,1){10}}
\put( 2,11.8){\bf 1}
\put(12,11.8){\bf 2}
\put(22,11.8){\bf 3}
\put(2, 1.8){\bf 4}

\end{picture}
\ *\ 4 \ =\
\left\{\rule{0pt}{16pt}\
\setlength{\unitlength}{1.3pt}%
\begin{picture}(70,20)(0,8)
\thicklines
\put(0, 0){\line(1,0){10}}
\put(0,10){\line(1,0){70}}
\put(0,20){\line(1,0){70}}
\put( 0, 0){\line(0,1){20}}
\put(10, 0){\line(0,1){20}}
\put(20,10){\line(0,1){10}}
\put(30,10){\line(0,1){10}}
\put(40,10){\line(0,1){10}}
\put(50,10){\line(0,1){10}}
\put(60,10){\line(0,1){10}}
\put(70,10){\line(0,1){10}}

\put(0,0){\usebox{\Shading}}

\put(32,11.8){\bf 5}
\put(42,11.3){\bf 6}
\put(52,11.3){\bf 7}
\put(62,11.3){\bf 8}

\end{picture}
\,,\
\setlength{\unitlength}{1.3pt}%
\begin{picture}(60,20)(0,8)
\thicklines
\put(0, 0){\line(1,0){20}}
\put(0,10){\line(1,0){60}}
\put(0,20){\line(1,0){60}}
\put( 0, 0){\line(0,1){20}}
\put(10, 0){\line(0,1){20}}
\put(20, 0){\line(0,1){20}}
\put(30,10){\line(0,1){10}}
\put(40,10){\line(0,1){10}}
\put(50,10){\line(0,1){10}}
\put(60,10){\line(0,1){10}}

\put(0,0){\usebox{\Shading}}

\put(12, 1.6){\bf 5}
\put(32,11.8){\bf 6}
\put(42,11.3){\bf 7}
\put(52,11.3){\bf 8}
\end{picture}
\,,\
\setlength{\unitlength}{1.3pt}%
\begin{picture}(50,20)(0,8)
\thicklines
\put(0, 0){\line(1,0){30}}
\put(0,10){\line(1,0){50}}
\put(0,20){\line(1,0){50}}
\put( 0, 0){\line(0,1){20}}
\put(10, 0){\line(0,1){20}}
\put(20, 0){\line(0,1){20}}
\put(30, 0){\line(0,1){20}}
\put(40,10){\line(0,1){10}}
\put(50,10){\line(0,1){10}}

\put(0,0){\usebox{\Shading}}

\put(12, 1.6){\bf 5}
\put(22, 1.3){\bf 6}
\put(32,11.8){\bf 7}
\put(42,11.3){\bf 8}

\end{picture}\
\right\}.
$$

Let $T(\alpha)\in T*\alpha$ be the tableau whose last
$\alpha$ entries lie in the first row.
Let $T^{+\alpha}\in T*\alpha$ be the tableau whose last $\alpha$ entries
lie in the second row, if that is possible.
Write $T^+$ for $T^{+1}$ and
define $T^{+\alpha}(\beta)$ to be $(T^{+\alpha})(\beta)$.
Then
$$
T*\alpha = \{T(\alpha),\, T^{+}(\alpha-1),\ldots,T^{+\alpha}\}.
$$

Given a set ${\cal T}$ of tableaux, define ${\cal T}*\alpha$ to be
$\bigcup_{T\in {\cal T}} T*\alpha$,
a disjoint union.
Similarly define ${\cal T}(\alpha)$, ${\cal T}^{+\alpha}$,
and ${\cal T}^{+\alpha}(\beta)$.
For $0\leq s\leq \alpha$, set
${\cal T}_{s,\alpha} = {\cal T}(\alpha) \bigcup {\cal T}^{+}(\alpha -1)
\bigcup\cdots\bigcup {\cal T}^{+s}(\alpha - s)$.
It follows that
${\cal T}_{s,\alpha} = {\cal T}_{s-1,\alpha} \bigcup
 \left( {\cal T}^{+s} (\alpha-s)\right)$
and ${\cal T}*\alpha = {\cal T}_{\alpha,\alpha}$.
Finally, for positive integers $\alpha_1,\ldots,\alpha_m$, define
$\alpha_1*\cdots*\alpha_m$ to be
$\emptyset*\alpha_1*\cdots*\alpha_m$, where $\emptyset$
is the unique tableau of shape $(0,0)$.

\subsection{Arrangements}\label{sec:arrangements}
An {\em arrangement} $\calf$ is a collection of $2n-3$ hyperplanes
$H_2,\ldots,H_{2n-2}$ and a complete flag
$\Fdot$  in $\Pn$,
which satisfy some additional assumtions.
For a tableau $T$ of degree at most $2n-2$, let $H_T$ be
$\bigcap \{H_j\,|\, j \in\mbox{ the second row of $T$}\}$.
For  $p=2,3,\ldots,2n-2$, the linear spaces in an arrangement
are required to
satisfy:
\begin{enumerate}
\item  $H_p \cap  F_{\hfl{p}} =  F_{\hfl{p}+1}$,
\item  For tableaux $S, T$ of degrees at most $p-1$,
$H_T\cap H_S\subset  H_p \Rightarrow H_T\cap H_S=F_{\hfl{p}+1}$.
\end{enumerate}
\noindent
The only obstruction to constructing arrangements over $k$ is that
$k$ have enough elements to find $H_p$ satisfying condition 2.
In particular, there are arrangements for any infinite field.
In \S\ref{thm:arrangement_finite_field}, we give an
estimate over which finite
fields it is possible to construct arrangements.

%
%

Let $\calf$ be an arrangement and  $T$ a tableau of shape
$(\alpha,\beta)$ with $\alpha + \beta \leq 2n-2$.
Define $\Omega(T) = \Omega(F_{\alpha+1}, H_T)$.
If $\alpha\geq n$, then $F_{\alpha+1} = \emptyset$, so
$\Omega(T) = \emptyset$.

\subsection{Lemma}\label{lemma:arrangements}
{\em Let $\calf$ be an arrangement,
$S,T$ tableaux of degree $l\leq 2n{-}2$ and suppose $T$ has shape
$(\alpha,\beta)$ with $\alpha\leq n-1$.
Then
\begin{enumerate}
\item $F_{\alpha+1} \subset F_{\hfl{l+1}+1} \subset H_T$
and $F_{\alpha+1}\subset H_T$ is a partial flag of type $(\alpha,\beta)$.
\item  If $\beta>0$, then $H_T\neq F_\beta$.
\item If  either $\Omega(T) =\Omega(S)$ or
$H_T=H_S$, then $T=S$.
\end{enumerate}
}
\medskip

\noindent{\bf Proof:}
Since $\alpha\geq \hfl{l{+}1}$ and for $p\leq l$,
$F_{\hfl{l{+}1}}\subset H_p$, we see
$F_{\alpha+1}\subset F_{\hfl{l{+}1}{+}1}\subset H_T$,
by the definition of arrangement.
We show the codimension of $H_T$ is $\beta$
by induction on $\beta$.
When $\beta = 0$ or $1$, this is clear.
Removing, if necessary, extra entries from the first row of $T$,
we may assume $l$ is in the second row of $T$.
Let $S$ be the tableau obtained from $T$ by removing $l$.
By induction, $H_S$ has codimension $\beta-1$.
By the definition of arrangement, $H_S \not\subset H_l$.
Noting $H_T = H_S\cap H_l$ proves (1).

For (2), if $H_T=F_\beta$, then $F_\beta\subset H_{p}$, where $p$ is
the $\beta$th entry in the second row of $T$.
Thus $F_\beta \subset F_{\hfl{p}+1}$ so
$\beta\geq\hfl{p}{+}1$.
But this contradicts $p \geq 2\beta$.

For (3), since $\Omega(T)\neq \emptyset$,
$H_T$ is the union of all lines in $\Omega(T)$.
Thus either assumption implies $H_T=H_S$.
Suppose $|S|+|T|$ is minimal subject to $S\neq T$
but $H_S = H_T$.
Let $s = |S|$ and $t=|T|$, then $s$ and $t$ are necessarily in
the second rows of $S$ and $T$, respectively.
If $s\neq t$, suppose $s>t$.
Then $H_T = H_S\cap H_s$, which contradicts
$H_T \not\subset H_s$.

Thus $s=t$.
Let $S'$ and $T'$ be the tableaux obtained by removing $s$ from
each of $S$ and $T$.
Then $H_{S'}\neq H_{T'}$, but $H_{S'}\cap H_{T'} = H_S \subset H_s$.
Thus $H_S = F_{\hfl{s}+1}$ a contradiction, as
the codimension of $H_S$ is at most $\hfl{s}$.
\QED

\subsection{The Cycles $\Omega({\cal T})$ and
$\Omega({\cal T}_{s,\alpha};L)$}
Let $\calf$ be an  arrangement.
For a set ${\cal T}$ of tableaux with common degree $l\leq 2n{-}2$,
define $\Omega({\cal T})$ to be the sum  over $T\in {\cal T}$ of
the Schubert varieties  $\Omega(T)$.
By  Lemma~\ref{lemma:arrangements}(3), these   are
distinct, so  $\Omega({\cal T})$ is a multiplicity
free cycle.

Let $0\leq s\leq \alpha$ be integers and $L$
a subspace of codimension $\alpha-s+1$
which meets the subspaces $F_{\hfl{l+s}+1},\ldots,F_{l+1}$
properly with $L\cap F_{l+1} = F_{l+\alpha-s+2}$ and
if $|T| = l$, then $L$ meets $H_{T^{+(s-1)}}$ and $H_{T^{+s}}$
properly.
When this occurs, we shall say that $L$
{\em meets ${\cal F}_{l,s}$ properly}.
Let $T\in {\cal T}_{s-1,\alpha}$.
Then $T$ is a tableau whose first row has length
$b\geq\hfl{l+s}+1+\alpha-s$ and no entry in its second row
exceeds $l+s-1$.
Define
$$
\Omega(T;L) = \left\{\begin{array}{ll}
\Omega(F_{b+1},H_T) & \mbox{ if } b\geq l+1\\
\Omega(F_{b-\alpha+s}\cap L, H_T) & \mbox{ otherwise }
\end{array}
\right. .
$$
Then we define  $\Omega({\cal T}_{s-1,\alpha};L)$
to be the sum over $T\in {\cal T}_{s-1,\alpha}$
of the varieties $\Omega(T;L)$.
By  Lemma~\ref{lemma:arrangements}(3), these   are
distinct, so $\Omega({\cal T}_{s-1,\alpha};L)$ is a multiplicity
free cycle.

\subsection{The sets $U_{i,s}$ and ${\cal G}_{i,s}$.}%
\label{sec:sets_of_cycles}
Let $\alpha_1,\ldots,\alpha_m$ be positive integers and
$\calf$ an arrangement.
Fix $1\leq i\leq m$.
Let ${\cal T}=\alpha_1*\cdots*\alpha_{i-1}$ and
$l$ be the common degree of  tableaux in
${\cal T}$.

Let $\Sigma_{i,0}\subset \G_1\Pn\times\prod_{j=i}^m\G^{\alpha_j+1}\Pn$
be the subvariety whose fibre over $(K_i,\ldots,K_m)$ is
$$
\Omega({\cal T})
\bigcap \Omega_{K_i}\bigcap\cdots\bigcap\Omega_{K_m}.
$$
Let $U_{i,0} \subset \prod_{j=i}^m\G^{\alpha_j+1}\Pn$
be those points for which this intersection is generically
transverse.
Let ${\cal G}_{i,0}\subset \Chow \G_1\Pn$ be
the fundamental cycles of fibres of $\Sigma_{i,0}$ over $U_{i,0}$.

For $1\leq s\leq \alpha_i$, let
$\Sigma_{i,s}\subset \G_1\Pn\times
\G^{\alpha_i-s+1}(\Pn/F_{l+\alpha_i-s+2})
\times \prod_{j=i+1}^m\G^{\alpha_j+1}\Pn$
be the subvariety where $L$ meets ${\cal F}_{l,s}$ properly and
whose fibre over $(L,K_{i+1},\ldots,K_m)$
is
$$
\left[\Omega({\cal T}_{s-1,\alpha_i};L)+
\Omega({\cal T}^{+s})\bigcap\Omega_L\right]
\bigcap \Omega_{K_{i+1}}
\bigcap\cdots\bigcap\Omega_{K_m}.
$$
Define $U_{i,s}$
and ${\cal G}_{i,s}$ analogously to $U_{i,0}$ and ${\cal G}_{i,0}$.

Set ${\cal G}_{m+1,0}$ to be the singleton
$\{ \Omega(\alpha_1*\cdots*\alpha_m)\}$.

When $k = \R$,
let ${\cal G}_{i,s;\Rsc}$ be the fundamental cycles of fibres of
$\Sigma_{i,s}$ over $U_{i,s}(\R)$.

\section{Main Results}

In \S4, we prove
\medskip

\noindent{\bf Theorem D.}\ {\em
Let $\alpha_1,\ldots,\alpha_m$ be positive
integers and $\calf$ any arrangement.
Then
\begin{enumerate}
\item For all $1\leq i\leq m$, and
$0\leq s\leq \alpha_i$, $U_{i,s}$ is a dense
open subset of the corresponding product of Grassmannians.
\item For all $1\leq i\leq m$, $0\leq s\leq \alpha_i$,
${\cal G}_{i,s}$ is a unirational subset of $\Chow\Gn$.
When $k=\R$, ${\cal G}_{i,s;\Rsc}$ is a real unirational
subset of $\Chow\Gn(\R)$.
\end{enumerate}
}
\medskip

Let $U$ be an open subset of $\PI$ and $\phi:U \rightarrow U_{i,s}$.
Then $\phi^*\Sigma_{i,s}$ has equidimensional  fibres over $U$.
As in \S\ref{sec:Chow}, the association of a point $u$ of $U$ to the
fundamental cycle of the fibre $(\phi^*\Sigma_{i,s})_u$ is
an algebraic morphism,
which we denote $\phi_* : U\rightarrow \Chow \Gn$.
Let $\phi_*$ be the unique extension of $\phi_*$ to
$\PI$, as well.
In \S5, we prove:
\medskip

\noindent{\bf Theorem E.}\ {\em
Let $\alpha_1,\ldots,\alpha_m$ be positive
integers and $\calf$ any arrangement.
Then
\begin{enumerate}
\item For all $1\leq i\leq m$, if $X$ is a closed point of
${\cal G}_{i+1,0}$, then there is an open subset $U$ of
$\PI- \{0\}$ and a map $\phi: U \rightarrow U_{i,\alpha_i}$
such that $X= \phi_*(0)$.
\item For all $1\leq i\leq m$ and $0\leq s\leq \alpha_i$,
if $X$ is a closed point of
${\cal G}_{i,s+1}$, then there is an open subset $U$ of
$\,\PI- \{0\}$ and a map $\phi: U \rightarrow U_{i,s}$
such that $X= \phi_*(0)$.
\end{enumerate}
}\medskip

\subsection{Proof of Theorem A}
In the situation of Theorem E, let ${\cal T} =\alpha_1*\cdots*\alpha_m$.
Then ${\cal G}_{1,0}$ is the set ${\cal G}$
and ${\cal G}_{m+1,0}$ is the cycle $\Omega({\cal T})$.
By Theorem~D part 2, any two points of ${\cal G}_{i,s}$ are connected
by a chain of rational curves, each lying within the closure of
${\cal G}_{i,s}$.
Downward induction in the lexicographic order on pairs $(i,s)$
gives a chain of rational curves between
$\Omega({\cal T})$ and a cycle $X\in {\cal G}$.
Thus Theorem E implies Theorem A
when $\beta_i = 0$ for $1\leq i \leq m$.

Suppose $\lambda^i = (\alpha_i+\beta_i,\beta_i)$ for $1\leq i \leq m$.
Let $M_0\subset \Pn$ have codimension
$\beta=\beta_1+\cdots+\beta_m$ and $\calf$ be any arrangement
in $M_0$.
Define $U_{i,s}$ and ${\cal G}_{i,s}$ as in \S\ref{sec:sets_of_cycles},
with $M_0$ replacing $\Pn$.

Let $X\in {\cal G}$,  so $X$ is a generically transverse intersection
$$
\Omega(K_1,M_1)\bigcap\cdots\bigcap\Omega(K_m,M_m),
$$
where $K_i\subset M_i$ has type $\lambda^i$ for $1\leq i\leq m$.
Set $M = M_1\cap\cdots\cap M_m$.
Iteration of Lemma~\ref{lemma:one} shows that $M$ has codimension
$\beta$ and $L_i = M \cap K_i$ has codimension $\alpha_i +1$ in $M$.
Thus
$$
X = \Omega_{L_1}\bigcap\cdots\bigcap \Omega_{L_m}
$$
is a generically transverse intersection in $\GI M$.

Let $\gamma$ be any automorphism of $\Pn$ with $\gamma M = M_0$
and  $\Gamma$ a one parameter subgroup containing $\gamma$.
The orbit $\Gamma \cdot X$ is a rational curve
(or a point) in $\Chow \Gn$  containing $\gamma(X)$.
Since $\gamma(X)$ is in ${\cal G}_{i,0}$,
previous arguments have shown there exists
a chain of rational curves between $\gamma(X)$ and $\Omega({\cal T})$,
each contained within the closure of ${\cal G}_{i,0}$.
\QED

\subsection{The Schubert Calculus.}
Let $\lambda^1,\ldots,\lambda^m$ be partitions with
$\lambda^i = (\alpha_i+\beta_i,\beta_i)$.
Set $\beta = \beta_1+\cdots+\beta_m$ and
$s_i = \alpha_i+\cdots+\alpha_{i-1}$.
For a partition $\lambda$ with $|\lambda| = s_m$, define
${\cal C}^\lambda = {\cal C}^\lambda(\alpha_1,\ldots,\alpha_m)$
to be those tableaux of shape $\lambda$
such that for $1\leq i\leq m$ the integers
$s_i+1,\ldots, s_i+\alpha_i$ occur in distinct columns
increasing from left to right.
Then $\alpha_1*\cdots*\alpha_m =
\bigcup_{|\lambda| = s_m} {\cal C}^\lambda$.
Let $c^\lambda_{\alpha_1,\ldots,\alpha_m} = \#{\cal C}^\lambda$.

Interpreting Theorems A and E in terms of products in the \medskip
Chow ring of $\Gn$, we have:

\noindent{\bf Theorem B$'$.}\
{\em
\begin{enumerate}
\item ${\displaystyle
\sigma_{\alpha_1}\cdots\sigma_{\alpha_m} =
\sum_{|\lambda| = s_m} c^\lambda_{\alpha_1,\ldots,\alpha_m} \sigma_\lambda}$.
\item ${\displaystyle
\sigma_{\lambda^1}\cdots\sigma_{\lambda^m} =
\sum_{\stackrel{\mbox{\scriptsize $\lambda = (\beta+a,\beta+b)$}}{a+b = s_m}}
c^{(a,b)}_{\alpha_1,\ldots,\alpha_m} \sigma_\lambda}$.
\end{enumerate}
}
\medskip

In particular, if $\alpha_1+\cdots +\alpha_m = 2n-2$, then
the only non-zero term on the right hand side of (1) is
$c^{(n-1,n-1)}_{\alpha_1,\ldots,\alpha_m}\sigma_{(n-1,n-1)}$, or
$c^{(n-1,n-1)}_{\alpha_1,\ldots,\alpha_m}$ times the class of a point (line).
Thus

\subsection{ Corollary}
{\em The number of lines meeting  general
$(n-\alpha_i - 1)$-planes for $1\leq i\leq m$ is equal to the number of
tableaux of shape $(n-1,n-1)$ such that for $1\leq i\leq m$, the integers
$s_i+1,\ldots, s_i+\alpha_i$ occur in distinct columns increasing from
left to right.
}
\medskip

This number, $c^{(n-1,n-1)}_{\alpha_1,\ldots,\alpha_m}$,
is also known as a Kostka number~\cite{Sagan}.

\subsection{Enumerative Geometry of the Real Grassmannian.}
Let ${\cal G}_{\Rsc}$ consist of the fundemental cycles
of generically transverse intersections of Schubert varieties
of types $\lambda^1,\ldots,\lambda^m$ defined by real flags.
\medskip

\noindent{\bf Theorem C$'$.} \ {\em
Let $\lambda^1,\ldots,\lambda^m$ be
partitions, suppose $\lambda^i = (\alpha_i+\beta_i,\beta_i)$,
and set $\beta=\sum_{i=1}^m \beta_i$.
Let $M \subset \Pn$ be a real $(n-\beta)$-plane and
$\calf$ an  arrangement in $M$.
\begin{enumerate}
\item
$\Omega(\alpha_1*\cdots*\alpha_m)$ is in the closure of
$\,{\cal G}_{\Rsc}$.
\item
If $\,|\lambda^1|+\cdots+|\lambda^m|=2n-2$,
then there is a nonempty classically open
subset  in the product of $m$ real flag manifolds whose
corresponding Schubert varieties meet transversally,
with all points of intersection real.
\end{enumerate}
}
\medskip

\noindent{\bf Proof:}
Suppose that $\lambda^i = (\alpha_i+\beta_i,\beta_i)$.
Define $U_{i,s}$ and ${\cal G}_{i,s}$ as in \S\ref{sec:sets_of_cycles}
for the arrangement ${\cal F}$ in $M$ and the integers
$\alpha_1,\ldots,\alpha_m$.
Arguing as in the proof of Theorem~A
shows ${\cal G}_{1,0;\Rsc}\subset {\cal G}_{\Rsc}$.
Restricting to the real points of the varieties in Theorem E
shows  ${\cal G}_{i,s;\Rsc}\subset \overline{{\cal G}_{1,0;\Rsc}}$.
The case $(i,s) = (m+1,0)$ is part 1.

For 2, let $d = c^{(n-1,n-1)}_{\alpha_1,\ldots,\alpha_m}$.
Then $\Omega(\alpha_1*\cdots*\alpha_m)$
is $d$ distinct real lines.
Hence ${\cal G}_{i,s} \subset S^d\Gn$, the Chow variety of
effective degree $d$ zero cycles on $\Gn$.
The real points $S^d\Gn(\R)$ of $S^d\Gn$ are
effective degree $d$ zero cycles
stable under complex conjugation.
The dense subset of $S^d\Gn(\R)$ of
multiplicity free cycles
has a component ${\cal M}$ parameterizing sets of $d$ distinct real lines
and $\Omega(\alpha_1*\cdots*\alpha_m)\in {\cal M}$.
By part 1,
$\Omega(\alpha_1*\cdots*\alpha_m)\in
\overline{{\cal G}_{\Rsc}}$, which shows
${\cal G}_{\Rsc}\bigcap{\cal M}\neq \emptyset$, a restatement
of 2.
\QED

\section{Generically Transverse  Intersections}

\subsection{Lemma.}\label{lemma:gic}
{\em
Let $\lambda^1,\ldots,\lambda^m$ be partitions.
Then the set $U$ of partial flags
$K_1\subset M_1,\,\ldots,$
$\,K_m\subset M_m$ for which the intersection
$$
\Omega(K_1,M_1)\bigcap\cdots\bigcap\Omega(K_m,M_m)
$$
is generically transverse is a dense open subset of
$\,\prod_{i=1}^m \Fl(\lambda^i)$.
}
\medskip

\noindent{\bf Proof:}
For $1\leq i\leq m$, let $K_i\subset M_i$ be a partial flag of
type $\lambda^i= (\alpha_i+\beta_i,\beta_i)$ and suppose the
corresponding Schubert varieties meet generically transversely.
By Lemma~\ref{lemma:one},
$$
\Omega(K_1,M_1)\bigcap\cdots\bigcap\Omega(K_m,M_m)
\ =\ \G_1M\bigcap \Omega_{L_1}\bigcap\cdots\bigcap\Omega_{L_m},
$$
where $M = M_1\cap\cdots\cap M_m$, $K_i = L_i\cap M_i$
where $L_i$ meets $M_i$ properly, and  $M$ has codimension
$\beta = \beta_1+\cdots+\beta_m$.
\medskip

Fix a codimension $\beta$ subspace $M$ of $\Pn$.
As $U$ is stable under the diagonal action of $Gl_{n+1}$, it is the
union of the translates of $V = U\cap X$, where
$X$ consists of those $m$-tuples of
flags with $M \subset M_i, 1\leq i\leq m$.
Moreover, $U$ is open if and only if $V$ is open in $X$.

Let $Y\subset X$ be those flags where
$M = M_1\cap\cdots\cap M_m$ and $K_i$ meets $M$ properly.
The product of maps defined by $(K_i,M_i) \mapsto K_i\cap M=L_i$
exhibits $Y$ as a fibre bundle over the product
$\prod_{i=1}^m \G^{\alpha_i+1}M$, and $V$ is the inverse image of the set
$W$ consisting of those $(L_1,\ldots,L_m)$ for which
$\Omega_{L_1}\bigcap\cdots\bigcap\Omega_{L_m}$ is generically transverse.
Thus we may assume $\beta_i = 0$.
\medskip

Let $\Sigma\subset(\Pn)^m\times\G_1\Pn\times\prod_{i=1}^m\G^{\alpha_i+1}\Pn$
consist of those
$(p_1,\ldots,p_m,\ell,L_1,\ldots,L_m)$ such that
$p_i\in \ell \cap L_i$ for $1\leq i\leq m$.
The projection of $\Sigma$  to $(\Pn)^m\times \G_1\Pn$
exhibits  $\Sigma$ as a  fibre bundle with fibre
$\prod_{i=1}^m\G^{\alpha_i+1}(\Pn/p_i)$
and  image  those $(p_1,\ldots,p_m,\ell)$ with each
$p_i\in \ell$.
This image has dimension $m + 2n-2$.
Thus $\Sigma$ is irreducible of dimension
$$
m+2n-2+\sum_{i=1}^m(n-\alpha_i-1)(\alpha_i+1)\ = \
2n-2 -\sum_{i=1}^m \alpha_i+ \sum_{i=1}^m (n-\alpha_i)(\alpha_i+1).
$$

The image of $\Sigma$ in $\prod_{i=1}^m\G^{\alpha_i+1}\Pn$ consists of those
$(L_1,\ldots,L_m)$ whose corresponding Schubert varieties have
nonempty intersection.
This image is a proper subvariety if
$2n-2 < \sum_{i=1}^m \alpha_i$.
In this case, $U$ is the complement of this image.

Suppose $2n-2 \geq \sum_{i=1}^m \alpha_i$.
Let $W\subset \Sigma$ consist of those points where
$\Omega_{L_1},\ldots,\Omega_{L_m}$ meet transversally at $\ell$.
By Lemma~\ref{lemma:one}, $W$ consists of those points such that
\begin{enumerate}
\item $\ell \not\subset L_i$ for $1\leq i\leq m$,
thus $p_i = \ell\cap L_i$ and $\ell$ is a smooth point of
$\Omega_{L_i}$.
\item The tangent spaces
$T_{\ell}\Omega_{L_i}$
meet transversally.
\end{enumerate}
Thus $W$ is an open subset of $\Sigma$.
We show $W\neq \emptyset$.
\smallskip

Fix $\ell \in \G_1\Pn$ and  distinct points
$p_1,\ldots,p_m$ of $\ell$.
Define $f : \Hom(\ell,\Pn/\ell)-\{0\} \rightarrow (\Pn/\ell)^m$
by $\phi \mapsto (\phi(p_1),\ldots,\phi(p_m))$.
Let $G\subset Gl_{n+1}$ fix $\ell$ pointwise.
Then $G^m$ acts transitively on $(\Pn/\ell)^m$.
Choose $(L_1,\ldots,L_m)\in
\prod_{i=1}^m\G^{\alpha_i+1}\Pn$
with $p_i = \ell \cap L_i$.

By Theorem 2 (i) of~\cite{Kleiman}, there is a dense open subset $V$ of
$G^m$ consisting of those {\boldmath $g$} such that
either
$f^{-1}(\mbox{\boldmath $g$}
((L_1+\ell)/\ell \times\cdots\times (L_m+\ell)/\ell))$
is empty or its codimension equals that of
$(L_1+\ell)/\ell \times\cdots\times (L_m+\ell)/\ell$ in $\Pn$,
which is $\sum_{i=1}^m \alpha_i$.

Let {\boldmath $g$} $ = (g_1,\ldots,g_m)\in V$ and set
$L_i' = g_i L_i$.
Then
$f^{-1}((L_1'+\ell)/\ell \times\cdots\times
(L_m'+\ell)/\ell)\cup \{0\}$
is the intersection of the tangent spaces
$T_\ell\Omega_{L_i'}$ for $1\leq i\leq m$.

Since $\alpha_i$ is the codimension of $T_\ell\Omega_{K_i'}$ in
$T_{\ell}\Gn$ and
$\sum_{i=1}^m \alpha_i \leq 2n-2$,
we see that
$\Omega_{L_1'},\ldots,\Omega_{L_m'}$ meet transversally at $\ell$.
Thus $W\neq \emptyset$.
\medskip

Let $Z = \Sigma - W$ and $\pi$ be the projection
$\Sigma\rightarrow \prod_{i=1}^m \G^{\alpha_i+1}\Pn$.
The desired set $U$
consists of those $(L_1,\ldots,L_m)$ with
$\dim(\pi^{-1}(L_1,\ldots,L_m) \bigcap Z)< 2n-2-\sum_{i=1}^m \alpha_i$.
$U$ is open and non-empty, for otherwise $\dim Z = \dim \Sigma$, which
implies $Z = \Sigma$ and contradicts $W\neq\emptyset$.
\QED

\subsection{Lemma.}\label{lemma:good_dimension}
{\em Let $d, \alpha_1,\ldots,\alpha_m$ be positive integers
and $Z$ a subscheme of $\Gn$ with $\dim(Z)<d$.
Then the set $W \subset \prod_{i=1}^m\G^{\alpha_i+1}\Pn$
consisting of those $(K_1,\ldots,K_m)$
for which
$\dim(Z\bigcap \Omega_{K_1}\bigcap\cdots\bigcap\Omega_{K_m})<
d-\sum_{i=1}^m\alpha_i$ is open and dense.
}
\medskip

\noindent{\bf Proof:}
Let $\Sigma\subset Z\times \prod_{i=1}^m\G^{\alpha_i+1}\Pn$
be the subscheme whose fibre over   $(K_1,\ldots,K_m)$
is
$Z\bigcap \Omega_{K_1}\bigcap\cdots\bigcap\Omega_{K_m}$.
By the upper semicontinuity of fibre dimension,
$W$ is open.
If $W$ were empty, then all fibres of the projection to
$\prod_{i=1}^m\G^{\alpha_i+1}\Pn$
would have dimension at least
$d-\sum_{i=1}^m\alpha_i$ and so
$\dim\Sigma \geq d-\sum_{i=1}^m\alpha_i +\sum_{i=1}^m
(n-\alpha_i)(\alpha_i+1)$.

Projecting  to $Z$ exhibits $\Sigma$ as a fibre bundle with fibre
over a point $\ell$ of $Z$ equal to $X_1(\ell)\times\cdots \times
X_m(\ell)$, where
$X_i(\ell)\subset \G^{\alpha_i+1}\Pn$ is the set of those
$K_i\subset\G^{\alpha_i+1}\Pn$ which meet $\ell$,
which has codimension $\alpha_i$.
Thus $\Sigma$ has dimension
$$
\dim Z  - \sum_{i=1}^m \alpha_i + \sum_{i=1}^m
(n-\alpha_i)(\alpha_i+1).
$$
Since $d> \dim Z$, $W$ must be non-empty.
\QED

\subsection{Proof of Theorem D, part 1}\label{sec:Proof_D}
We show that for each $1\leq i\leq m$ and
$0\leq s\leq \alpha_i$, the sets $U_{i,s}$ are open dense
subsets of the corresponding products of Grassmannians.

Let ${\cal T} = \alpha_1*\cdots*\alpha_{i-1}$.
Recall that $U_{i,0}$
consists of those $(K_i,\ldots,K_m)$ such that
the intersection
$$
\Omega({\cal T})
\bigcap \Omega_{K_i}\bigcap\cdots\bigcap\Omega_{K_m}
$$
is generically transverse.
Such a cycle has dimension $d = 2n-2 -\sum_{i=1}^m \alpha_i$.

Let $Z$ be  the singular locus
of $\Omega({\cal T})$.
The above intersection is generically transverse
if $\dim (Z \bigcap  \Omega_{K_i}\bigcap\cdots\bigcap\Omega_{K_m}) <d$
and if, for every component $\Omega(T)$ of $\Omega({\cal T})$,
the intersection
$\Omega(T)\bigcap  \Omega_{K_i}\bigcap\cdots\bigcap\Omega_{K_m}$
is generically transverse.

By Lemma~\ref{lemma:arrangements},
$\Omega(T) = \Omega(S)\neq \emptyset$ implies that
$T=S$, thus $Z$ is a union of intersections of components and
the singular loci of components, and hence
$\dim(Z) < \dim(\Omega({\cal T}))$.
By Lemma~\ref{lemma:good_dimension},
there is an open subset $W$ of $\prod_{j=i}^m \G^{\alpha_j+1}\Pn$
consisting those $(K_1,\ldots,K_m)$ for which
$\dim(Z\bigcap \Omega_{K_i}\bigcap\cdots\bigcap\Omega_{K_m})<d$.

For  $T\in {\cal T}$, let $U_T\subset \prod_{j=i}^m\G^{\alpha_j+1}\Pn$
be those $(K_i,\ldots,K_m)$ where the intersection
$$
\Omega(T)\bigcap \Omega_{K_i}\bigcap\cdots\bigcap\Omega_{K_m}
$$
is generically transverse.
It suffices to show that for each $T\in {\cal T}$,
$U_T$ is a dense open subset of $\prod_{j=i}^m\G^{\alpha_j+1}\Pn$,
since $U_{i,0} =
W\cap\bigcap_{T\in {\cal T}} U_T$.

Suppose $T$ has shape $\lambda = (\alpha,\beta)$ and let
$V\subset \Fl(\lambda)\times\prod_{j=i}^m\G^{\alpha_j+1}\Pn$
be those flags for which the intersection
$$
\Omega(F,H)\bigcap \Omega_{K_i}\bigcap\cdots\bigcap\Omega_{K_m}
$$
is generically transverse.
By Lemma~\ref{lemma:gic}, $V$ is dense and open.
Note that
$$
\{F_{\alpha+1}\subset H_T\}\times U_T =
V\bigcap\left( \{F_{\alpha+1}\subset H_T\}\times
\prod_{j=i}^m\G^{\alpha_j+1}\Pn\right),
$$
so $U_T$ is open.
Since $\Fl(\lambda) = Gl_{n+1} \cdot \{F_{\alpha+1}\subset H_T\}$,
and $V$ is stable under the diagonal
action of $Gl_{n+1}$,
$V =  Gl_{n+1} \cdot (\{F_{\alpha+1}\subset H_T\}\times U_T)$.
Thus $U_T$ is non-empty.
\medskip

The case of $U_{i,s}$ for $s>0$ follows by similar arguments.
\QED

\subsection{Unirationality of ${\cal G}_{i,s}$}
\mbox{ }\medskip

\noindent{\bf Lemma.}\label{lemma:unirationality}
{\em
Let $X$ be a projective variety, $U$  a dense open subset of a
variety $Y$ which is covered by affine spaces, and suppose that
$\Sigma \subset X\times U$ is closed and the projection to $U$
has generically reduced fibres of pure dimension.
Then the set ${\cal G}\subset \Chow X$ of
fundamental cycles of the fibres of $\,\Sigma$ is unirational.
When  $k=\R$, let ${\cal G}_{\Rsc}$ be those cycles
lying over $U(\R)$, then  ${\cal G}_{\Rsc}$ is
real unirational.
}
\medskip

Part 2 of Theorem D is a consequence of this Lemma.

\noindent{\bf Proof:}
Let $\pi$ be the projection $\Sigma \rightarrow U$.
As in \S 2.1, let $U'\subset U$ be the open set where the map $\phi$
which associates a point of $U$ to the fundamental cycle of the fibre
of $\pi$ at $x$ is an algebraic morphism.
Then $\phi(U')\subset {\cal G}$.
We show ${\cal G}\subset \overline{\phi(U')}$.

Let $x\in U$ and
choose a map $f:{\bf A}^1 \rightarrow Y$ with $f(0) = x$ and
$f^{-1}(U')\neq \emptyset$.
This is possible as $Y$ has a covering by affine spaces.
Then $f^{-1}\Sigma \rightarrow f^{-1}(U)$ is  a family over a smooth curve
with generically reduced fibres of pure dimension.
The association of a point $u$ of $f^{-1}(U)$ to the
fundamental cycle of the fibre $\pi^{-1}(f(u))$ gives a map
$\psi: f^{-1}(U) \rightarrow \Chow X$ agreeing with
$\phi\circ f$ on $f^{-1}(U')$.
Thus the fundamental cycle of $\pi^{-1}(x)$
is in $\phi(U')$.

If $k = \R$, these maps show
$\phi(U'(\R))\subset{\cal G}_{\Rsc}
\subset \overline{\phi(U'(\R))}$,
thus ${\cal G}_{\Rsc}$ is real unirational.
\QED

\section{Construction of Explicit Rational Equivalences}

We use the following to  parameterize the explicit rational
equivalences we construct.

\subsection{Lemma.}\label{lemma:limits_are_good}
{\em
Let $\Fdot$ be a complete flag in $\Pn$.
Suppose $L_{\infty}$ is a hyperplane not containing $F_n$.
Then there exists a pencil of hyperplanes $L_t$, for
$t\in \Pp^1 = {\bf A}^1\bigcup \{\infty\}$,
such that if $t\neq 0$, then $L_t$ meets the subspaces
in $\Fdot$ properly,
and,  for each $i\leq n-1$, the family of codimension $i+1$ planes
induced by $L_t\bigcap F_i$, for $t\neq 0$ has fibre $F_{i+1}$ over $0$.
}
\medskip

\noindent{\bf Proof:}
Let $x_0,\ldots,x_n$ be coordinates for $\Pn$ such that
$L_{\infty}$ is given by $x_n=0$ and
$F_i$ by $x_0=\cdots=x_{i-1}=0$.
Let $e_0,\ldots,e_n$ be a basis for $\Pn$ dual to these coordinates
and define
$$
L_t = \Span{te_j + e_{j+1}\,|\, 0\leq j\leq n-1}.
$$
For $t\neq 0$,
$L_t \bigcap F_i = \Span{te_j + e_{j+1}\,|\, i\leq j\leq n-1}$
and so has codimension $i+1$.
The fibre of this family at $0$ is
$\Span{e_{j+1}\,|\, i\leq j\leq n-1} = F_{i+1}$.
\QED

In the situation of Lemma~\ref{lemma:limits_are_good}, we write
$\lim_{t\rightarrow 0} L_t\cap F_i = F_{i+1}$.
\medskip

For the remainder of this section, fix an arrangement $\calf$.
Set ${\cal T} = \alpha_1*\cdots*\alpha_{i-1}$ and let
$l$ be  the common degree of tableaux in
${\cal T}$.

\subsection{Proof of Theorem E, part 1.}\label{sec:proof_E_1}
Let $1\leq i\leq m$,
and suppose  $X_0$ is a  cycle in ${\cal G}_{i+1,0}$:
$$
X_0 =  \Omega({\cal T}*\alpha_i)
\bigcap\Omega_{K_{i+1}}\bigcap\cdots\bigcap\Omega_{K_m}.
$$

Let $L_\infty$ be any hyperplane which meets ${\cal F}_{l,\alpha_i}$ properly.
By Lemma~\ref{lemma:limits_are_good} applied to the flag
induced by $\Fdot$ in $\Pn/F_{l+2}$ and the hyperplane $L_\infty/F_{l+2}$,
there is a pencil  $L_t$ of hyperplanes such that if $t\neq 0$
and $i\leq l+1$, then $L_t$ meets $F_i$ properly and
$\lim_{t\rightarrow 0} L_t\cap F_i = F_{i+1}$.

Let ${\cal X}\subset \Pp^1\times \Gn$ be the subscheme whose fibre
at $t\neq 0$ is
$$
X_t = \left[ \Omega({\cal T}_{\alpha_i-1,\alpha_i};L_t) + \rule{0pt}{13pt}
\Omega({\cal T}^{+\alpha})\right] \bigcap
\Omega_{K_{i+1}}\bigcap\cdots\bigcap\Omega_{K_m}.
$$
Since $L_\infty$ meets ${\cal F}_{l,\alpha_i}$ properly,
the set $U'\subset \PI$ of $t$  where $L_t$ meets
${\cal F}_{l,\alpha_i}$ properly is open and dense.

We claim that  $X_0$ is the fibre of ${\cal X}$ over $0$.
In that case, let $U''\subset \Pp^1$ be the open subset of those
$t$ for which
$X_t$ is generically reduced.
Since $X_0$ is generically reduced,
$0 \in U''$ so $U'' \neq \emptyset$.
For  $t\in U'\cap U''$, the fibre
$X_t\in {\cal G}_{i-1,\alpha_i}$,
as $\Omega_{L_t} = \Gn$.
The restriction of ${\cal X}$ to $U\cup \{0\}$ gives a family
over a smooth curve with generically reduced equidimensional fibres.
Thus the association of a point of $U\cup \{0\}$ to its fibre
gives a map $\phi: U\cup \{0\}\rightarrow \Chow\Gn$ with
$\phi(U)\subset {\cal G}_{i-1,\alpha_i}$ and $\phi(0) = X_0$,
proving part 1.

For $T\in  {\cal T}_{\alpha_i-1,\alpha_i}$, let
${\cal X}_T\subset  \Pp^1\times \Gn$ be the subscheme whose fibre
at $t\neq 0$ is
$$
(X_T)_t = \Omega(T;L_t)
\bigcap\Omega_{K_{i+1}}\bigcap\cdots\bigcap\Omega_{K_m}.
$$
Since ${\cal X}
= \sum_{T\in  {\cal T}_{\alpha_i-1,\alpha_i}} {\cal X}_T \ +\
\Pp^1\times \Omega({\cal T}^{+\alpha_i})
\bigcap\Omega_{K_{i+1}}\bigcap\cdots\bigcap\Omega_{K_m}$,
to show that $X_0$ is the fibre of ${\cal X}$ over $0$
it suffices to show that for each $T\in {\cal T}_{\alpha_i-1,\alpha_i}$,
the fibre of ${\cal X}_T$ at 0 is
$\Omega(T)\bigcap\Omega_{K_{i+1}}\bigcap\cdots\bigcap\Omega_{K_m}$.

Let  $T\in {\cal T}_{\alpha_i-1,\alpha_i}$.
If the first row of $T$ has length exceeding $l+1$, then
$\Omega(T;L_t) = \Omega(T)$, so ${\cal X}_T$ is the constant family
$\Pp^1\times\Omega(T)\bigcap\Omega_{K_{i+1}}\bigcap\cdots\bigcap\Omega_{K_m}$.

Now suppose the first row of $T$ has length $b\leq l+1$.
Then, for $t\neq 0$,  $\Omega(T;L_t) = \Omega(F_b\cap L_t,H_T)$.
Since $\lim_{t\rightarrow 0}F_b\bigcap L_t = F_{b+1}$, we see that
$\Omega(T)$ is the fibre over 0 of the family over ${\bf P}^1$
whose fibre over $t\neq 0$ is $\Omega(T;L_t)$.

Since
$\Omega(T)\bigcap\Omega_{K_{i+1}}\bigcap\cdots\bigcap\Omega_{K_m}$
is generically transverse, there is an open subset $U_T\subset {\bf P}^1$
such that for $t\in U_T-\{0\}$,
$\Omega(T;L_t)\bigcap\Omega_{K_{i+1}}\bigcap\cdots\bigcap\Omega_{K_m}$
is generically transverse.
This shows that the fibre over 0
of ${\cal X}_T$ is
$\Omega(T)\bigcap\Omega_{K_{i+1}}\bigcap\cdots\bigcap\Omega_{K_m}$.
\medskip

\subsection{Proof of Theorem E, part 2}

Let $0\leq s\leq \alpha_i-1$
and suppose $X_0\in {\cal G}_{i,s+1}$
$$
X_0 =  \left[\Omega({\cal T}_{s,\alpha_i};N)+
\Omega({\cal T}^{+s+1})\bigcap\Omega_N
\right]\bigcap\Omega_{K_{i+1}}\bigcap\cdots\bigcap\Omega_{K_m}.
$$
Then $N$ has codimension $\alpha_i-s$ in $\Pn$
and meets ${\cal F}_{l,s+1}$ properly,
and the above intersection is generically transverse.
We make a useful calculation.
Let $L_0 = N\cap H_{l+s+1}$.
\medskip

\noindent{\bf Lemma. }
$$
\Omega({\cal T}^{+s}(\alpha_i-s);N)+\Omega({\cal T}^{+s+1})\bigcap\Omega_N
= \Omega({\cal T}^{+s})\bigcap \Omega_{L_0}.
$$
\medskip

Since ${\cal T}_{s,\alpha_i} = {\cal T}_{s-1,\alpha_i}
\bigcup {\cal T}^{+s}(\alpha_i-s)$, we see that
$$
X_0 =  \left[\Omega({\cal T}_{s-1,\alpha_i};N)+
\Omega({\cal T}^{+s})\bigcap \Omega_{L_0}
\right]\bigcap \Omega_{K_{i+1}}\bigcap\cdots\bigcap\Omega_{K_m}.
$$

\noindent{\bf Proof:}
Let $T\in {\cal T}^{+s}$
and suppose that $b$ is the length of the first row of $T$.
Then $\hfl{l+s+1}\leq b\leq l$.
The degree of $T$ is $l+s$, so
$L_0 = N\cap H_{l+s+1}$ meets $H_T$ properly,
because $H_T\cap H_{l+s+1}$ equals either $H_{T^+}$ or $F_{\hfl{l+s+1}+1}$,
each of which meets $N$ properly.

If $T$ is rectangular,  $b = \frac{l+s}{2}$ and
$\Omega(T) = \G_1 H_T$.
Thus,
$\Omega(T)\bigcap \Omega_{L_0} =
\Omega( H_T\cap L_0, H_T)$.
Since $L_0$ meets $H_T$ properly, this is generically
transverse, by Lemma~\ref{lemma:one}.

We calculate $ H_T\cap L_0$.
First note that
$ H_T\cap H_{l+s+1}=F_{\hfl{l+s+1}+1} = F_{\hfl{l+s}+1}$.
So
$L_0\cap H_T = N\cap H_{l+s+1}\cap H_T
= N\cap F_{\hfl{l+s}+1}$.
As $H_T = H_{T(\alpha_i-s)}$,
$$
\Omega(F_{b+1},H_T)\bigcap \Omega_{L_0} =
\Omega(N\cap F_{\hfl{l+s}+1},H_{T(\alpha_i-s)}) = \Omega(T(\alpha_i-s);N).
$$

Suppose $T$ is not rectangular.
Since $F_{b+1}\subset F_{\hfl{l+s+1}+1}\subset H_{l+s+1}$,
Lemma~\ref{lemma:component_calculation} implies
$$
\Omega(F_{b+1},H_T)\bigcap \Omega_{L_0} =
\Omega(F_{b+1}\cap N,H_T) + \Omega(F_{b+1},H_T\cap H_{l+s+1})
\bigcap\Omega_N.
$$
But this is
$\Omega(T(\alpha_i-s);N) + \Omega(T^+)\bigcap \Omega_N$.

Summing over $T\in {\cal T}^{+s}$ completes the proof.
\QED

Let $\Ndot$ be any complete flag in
$N/F_{l+\alpha_i-s+2}$ refining the images of
$$
F_{l+\alpha_i-s+1}\subset N\cap F_l\subset \cdots\subset
N\cap F_{\hfl{l+s}+1}\subset L_0.
$$
Let $L_\infty$ be any hyperplane of $N$ which meets ${\cal F}_{l,s}$
properly.
By Lemma~\ref{lemma:limits_are_good}
applied to $\Ndot$ in $N/F_{l+\alpha_i-s+2}$, there is a
pencil $L_t$ of hyperplanes
of $N$ each containing $F_{l+\alpha_i-s+2}$,
and  for $\hfl{l+s}+1\leq j\leq l$ and $t\neq 0$,
$L_t$ meets $N\cap F_j$ properly, with
$\lim_{t\rightarrow 0} L_t\cap N\cap  F_j
= N\cap F_{j+1}$.
Since $L_t\cap N\cap F_j = L_t\cap F_j$ for $j$ in this
range, $L_t$ meets $F_j$ properly.

Let ${\cal X}\subset \PI\times\Gn$ be the subscheme whose fibre
over $t\in \PI$ is
$$
X_t = \left[\Omega({\cal T}_{s-1,\alpha_i};L_t) +
\Omega({\cal T}^{+s})\bigcap \Omega_{L_t}\right] \bigcap
\Omega_{K_{i+1}}\bigcap\cdots\bigcap\Omega_{K_m}.
$$
Since $L_\infty$ meets ${\cal F}_{l,\alpha_i}$ properly,
The set $U'\subset \PI$ of $t$  where $L_t$ meets
${\cal F}_{l,\alpha_i}$ properly is open and dense.

We claim that $X_0$ is the fibre of ${\cal X}$ over $0$.
In that case, let $U''\subset \Pp^1$ be the open subset of those
$t$ for which
$X_t$ is generically reduced.
Since $X_0$ is generically reduced,
$0 \in U''$ so $U'' \neq \emptyset$.
For  $t\in U'\cap U''$, the fibre
$X_t\in {\cal G}_{i,s}$.
The restriction of ${\cal X}$ to $U\cup \{0\}$ gives a family
over a smooth curve with generically reduced equidimensional fibres.
Thus the association of a point of $U\cup \{0\}$ to its fibre
gives a map $\phi: U\cup \{0\}\rightarrow \Chow\Gn$ with
$\phi(U)\subset {\cal G}_{i,s}$ and $\phi(0) = X_0$,
proving part 1.

For $T\in {\cal T}_{s-1,\alpha_i}$
let ${\cal X}_T$ be the subscheme of $\PI\times \Gn$ whose fibre over
$t\neq 0$ is
$$
({\cal X}_T)_t =
\Omega(T;L_t)\bigcap \Omega_{K_{i+1}}\bigcap\cdots\bigcap\Omega_{K_m}.
$$
Arguing as at the end of \S\ref{sec:proof_E_1},
we may conclude that
$\Omega(T;N)\bigcap \Omega_{K_{i+1}}\bigcap\cdots\bigcap\Omega_{K_m}$
is the fibre of ${\cal X}_T$ at 0.

For $S\in {\cal T}^{+s}$, let ${\cal X}_S$ be the subscheme of
$\PI\times \Gn$ whose fibre over $t$ is
$$
({\cal X}_S)_t =
\Omega(S)\bigcap \Omega_{L_t} \bigcap
\Omega_{K_{i+1}}\bigcap\cdots\bigcap\Omega_{K_m}.
$$
Arguing as  at the end of \S\ref{sec:proof_E_1},
we may conclude that
$\Omega(S)\bigcap \Omega_{L_0} \bigcap
\Omega_{K_{i+1}}\bigcap\cdots\bigcap\Omega_{K_m}$
is the fibre of ${\cal X}_S$ at 0.

Since ${\cal X} = \sum_{T\in {\cal T}_{s-1,\alpha_i}}
{\cal X_T}  +\sum_{S\in {\cal T}^{+s}}{\cal X}_S$,
we conclude that  the fibre of ${\cal X}$ at 0 is $X_0$.
\QED

\section{An Algebra of Tableaux}

The Schubert classes, $\sigma_{\lambda}$, form an integral basis for
the Chow ring of any Grassmann variety.
Thus there exist integral constants
$c^{\lambda}_{\mu\,\nu}$ defined by the identity:
$$
\sigma_{\mu}\cdot\sigma_{\nu} =
\sum_\lambda c^{\lambda}_{\mu\,\nu}\sigma_{\lambda}.
$$
In 1934,
Littlewood and Richardson~\cite{Littlewood_Richardson} gave a conjectural
formula for these constants,
which was proven in 1978 by Thomas~\cite{Thomas_schensted_construction}.

Lascoux and Sch\"utzenberger~\cite{Lascoux_Schutzenberger_monoid_plactic}
constructed the ring of symmetric functions as a subalgebra of a
non-commutative
associative ring called the plactic algebra whose additive group is the
free abelian group $\Lambda$ with basis the set of Young tableaux.
For that, each tableau $T$ of shape $\lambda$
determines a monomial summand of the
Schur function, $s_{\lambda}$.

Evaluating $s_{\lambda}$ at  Chern roots of the dual to
the tautological bundle of the Grassmannian gives the Schubert class
$\sigma_\lambda$.
Non-symmetric monomials in these Chern roots are not defined,
so individual Young tableaux are not expected to appear in the
geometry of Grassmannians.
In this context, the crucial use we made of the Schubert varieties
$\Omega(T)$ is surprising.

A feature of our methods is the correspondence
between an iterative construction of the set
$\alpha_1*\cdots*\alpha_m$ and the rational
curves in the proof of Theorem~E.
This suggests an alternate non-commutative
associative product $\circ$ on $\Lambda$.
The resulting algebra has surjections
to the ring of symmetric functions and to Chow rings of
Grassmannians.

Additional combinatorial preliminaries for this section
may be found in~\cite{Sagan}.
Here, partitions $\lambda$, $\mu$, and $\nu$ may have any
number of rows.
Suppose $T$ and $U$ are,
respectively, a tableau of shape $\mu$ and a skew tableau
of shape $\lambda/\mu$.
Let $T\bigcup U$ be the tableau of shape $\lambda$ whose
first $|\mu|$ entries comprise $T$, and remaining entries
comprise $U$, with each increased by $|\mu|$.
For tableaux $S$ and $T$ where the shape of $S$ is $\lambda$, define
$$
S \circ  T = \sum S\bigcup U,
$$
the sum over all $\nu$ and all skew tableaux $U$
of shape $\nu/\lambda$  Knuth equivalent to $T$.
For example:
%
%
\begin{eqnarray*}
\setlength{\unitlength}{1.3pt}%
\begin{picture}(30,10)(0,12)
\thicklines
\put(0, 0){\line(1,0){10}}
\put(0,10){\line(1,0){30}}
\put(0,20){\line(1,0){30}}
\put( 0, 0){\line(0,1){20}}
\put(10, 0){\line(0,1){20}}
\put(20,10){\line(0,1){10}}
\put(30,10){\line(0,1){10}}
\put( 2,11.8){\bf 1}
\put(12,11.8){\bf 2}
\put(22,11.8){\bf 3}
\put( 2, 1.8){\bf 4}
\end{picture}
\ \circ\
\begin{picture}(40,10)(0,2)
\thicklines
\put(0, 0){\line(1,0){40}}
\put(0,10){\line(1,0){40}}
\put( 0, 0){\line(0,1){10}}
\put(10, 0){\line(0,1){10}}
\put(20, 0){\line(0,1){10}}
\put(30, 0){\line(0,1){10}}
\put(40, 0){\line(0,1){10}}
\put( 2, 1.8){\bf 1}
\put(12, 1.8){\bf 2}
\put(22, 1.8){\bf 3}
\put(32, 1.8){\bf 4}
\end{picture}
 &=&
\setlength{\unitlength}{1.3pt}%
\begin{picture}(70,10)(0,12)
\thicklines
\put(0, 0){\line(1,0){10}}
\put(0,10){\line(1,0){70}}
\put(0,20){\line(1,0){70}}
\put( 0, 0){\line(0,1){20}}
\put(10, 0){\line(0,1){20}}
\put(20,10){\line(0,1){10}}
\put(30,10){\line(0,1){10}}
\put(40,10){\line(0,1){10}}
\put(50,10){\line(0,1){10}}
\put(60,10){\line(0,1){10}}
\put(70,10){\line(0,1){10}}

\put(0,0){\usebox{\Shading}}
\put(32,11.8){\bf 5}
\put(42,11.3){\bf 6}
\put(52,11.3){\bf 7}
\put(62,11.3){\bf 8}

\end{picture}
\ +\
\setlength{\unitlength}{1.3pt}%
\begin{picture}(60,10)(0,12)
\thicklines
\put(0, 0){\line(1,0){20}}
\put(0,10){\line(1,0){60}}
\put(0,20){\line(1,0){60}}
\put( 0, 0){\line(0,1){20}}
\put(10, 0){\line(0,1){20}}
\put(20, 0){\line(0,1){20}}
\put(30,10){\line(0,1){10}}
\put(40,10){\line(0,1){10}}
\put(50,10){\line(0,1){10}}
\put(60,10){\line(0,1){10}}

\put(0,0){\usebox{\Shading}}
\put(12, 1.6){\bf 5}
\put(32,11.8){\bf 6}
\put(42,11.3){\bf 7}
\put(52,11.3){\bf 8}
\end{picture}
\ +\
\setlength{\unitlength}{1.3pt}%
\begin{picture}(50,10)(0,12)
\thicklines
\put(0, 0){\line(1,0){30}}
\put(0,10){\line(1,0){50}}
\put(0,20){\line(1,0){50}}
\put( 0, 0){\line(0,1){20}}
\put(10, 0){\line(0,1){20}}
\put(20, 0){\line(0,1){20}}
\put(30, 0){\line(0,1){20}}
\put(40,10){\line(0,1){10}}
\put(50,10){\line(0,1){10}}

\put(0,0){\usebox{\Shading}}

\put(12, 1.6){\bf 5}
\put(22, 1.3){\bf 6}
\put(32,11.8){\bf 7}
\put(42,11.3){\bf 8}
\end{picture}
\ + \raisebox{-12pt}{\rule{0pt}{5pt}}
\\
& &
\setlength{\unitlength}{1.3pt}%
\begin{picture}(70,20)(0,12)
\thicklines
\put(0,-10){\line(1,0){10}}
\put(0, 0){\line(1,0){10}}
\put(0,10){\line(1,0){60}}
\put(0,20){\line(1,0){60}}
\put( 0,-10){\line(0,1){30}}
\put(10,-10){\line(0,1){30}}
\put(20,10){\line(0,1){10}}
\put(30,10){\line(0,1){10}}
\put(40,10){\line(0,1){10}}
\put(50,10){\line(0,1){10}}
\put(60,10){\line(0,1){10}}
\put(0,0){\usebox{\Shading}}
\put( 2,-8.2){\bf 5}
\put(32,11.3){\bf 6}
\put(42,11.3){\bf 7}
\put(52,11.3){\bf 8}

\end{picture}
\ +\
\setlength{\unitlength}{1.3pt}%
\begin{picture}(60,20)(0,12)
\thicklines
\put(0,-10){\line(1,0){10}}
\put(0, 0){\line(1,0){20}}
\put(0,10){\line(1,0){50}}
\put(0,20){\line(1,0){50}}
\put( 0,-10){\line(0,1){30}}
\put(10,-10){\line(0,1){30}}
\put(20, 0){\line(0,1){20}}
\put(30,10){\line(0,1){10}}
\put(40,10){\line(0,1){10}}
\put(50,10){\line(0,1){10}}
\put(0,0){\usebox{\Shading}}
\put( 2,-8.2){\bf 5}
\put(12, 1.8){\bf 6}
\put(32,11.3){\bf 7}
\put(42,11.3){\bf 8}
\end{picture}
\ +\
 \setlength{\unitlength}{1.3pt}%
\begin{picture}(50,20)(0,12)
\thicklines
\put(0,-10){\line(1,0){10}}
\put(0, 0){\line(1,0){30}}
\put(0,10){\line(1,0){40}}
\put(0,20){\line(1,0){40}}
\put( 0,-10){\line(0,1){30}}
\put(10,-10){\line(0,1){30}}
\put(20, 0){\line(0,1){20}}
\put(30, 0){\line(0,1){20}}
\put(40,10){\line(0,1){10}}

\put(0,0){\usebox{\Shading}}

\put( 2,-8.2){\bf 5}
\put(12, 1.3){\bf 6}
\put(22, 1.8){\bf 7}
\put(32,11.3){\bf 8}
\end{picture}
\ .
\raisebox{-35pt}{\rule{0pt}{5pt}}
\end{eqnarray*}
This product is related to the composition $*$ of
\S\ref{sec:calculus_of_tableaux}:
Let $Y_\alpha$ be the unique standard
tableau of shape $(\alpha,0)$.
Then $T*\alpha$ consists of the summands of
$T\circ Y_{\alpha}$ with at most two rows.
\medskip

\noindent{\bf Theorem F.}\ {\em
The product  $\circ $ determines an associative non-commutative
$\bf Z$-algebra structure on $\Lambda$ with unit the empty tableau
$\emptyset$.
Moreover, $\circ$ is not the plactic product.
}
\medskip

\noindent{\bf Proof:}
In the plactic algebra, the product of two tableaux is always a third,
showing $\circ$ is not the plactic product.
For any tableau $T$, $\emptyset \circ T = T\circ \emptyset = T$.
Note
$$
\setlength{\unitlength}{1.3pt}%
\begin{picture}(10,10)(0,3)
\thicklines
\put(0, 0){\line(1,0){10}}
\put(0,10){\line(1,0){10}}
\put( 0, 0){\line(0,1){10}}
\put(10, 0){\line(0,1){10}}
\put( 2, 1.8){\bf 1}
\end{picture}
\ \circ\
\begin{picture}(20,8)(0,3)
\thicklines
\put(0, 0){\line(1,0){20}}
\put(0,10){\line(1,0){20}}
\put( 0, 0){\line(0,1){10}}
\put(10, 0){\line(0,1){10}}
\put(20, 0){\line(0,1){10}}
\put( 2, 1.8){\bf 1}
\put(12, 1.8){\bf 2}
\end{picture}
\ = \
\begin{picture}(30,10)(0,3)
\thicklines
\put(0, 0){\line(1,0){30}}
\put(0,10){\line(1,0){30}}
\put( 0, 0){\line(0,1){10}}
\put(10, 0){\line(0,1){10}}
\put(20, 0){\line(0,1){10}}
\put(30, 0){\line(0,1){10}}
\put( 2, 1.8){\bf 1}
\put(12, 1.8){\bf 2}
\put(22, 1.8){\bf 3}
\end{picture}
\ +\
\setlength{\unitlength}{1.3pt}%
\begin{picture}(20,20)(0,8)
\thicklines
\put(0, 0){\line(1,0){10}}
\put(0,10){\line(1,0){20}}
\put(0,20){\line(1,0){20}}
\put( 0, 0){\line(0,1){20}}
\put(10, 0){\line(0,1){20}}
\put(20,10){\line(0,1){10}}
\put( 2,11.8){\bf 1}
\put(12,11.8){\bf 3}
\put( 2, 1.8){\bf 2}
\end{picture}
\ \ \ \neq\ \ \
\begin{picture}(30,10)(0,3)
\thicklines
\put(0, 0){\line(1,0){30}}
\put(0,10){\line(1,0){30}}
\put( 0, 0){\line(0,1){10}}
\put(10, 0){\line(0,1){10}}
\put(20, 0){\line(0,1){10}}
\put(30, 0){\line(0,1){10}}
\put( 2, 1.8){\bf 1}
\put(12, 1.8){\bf 2}
\put(22, 1.8){\bf 3}
\end{picture}
\ +\
\setlength{\unitlength}{1.3pt}%
\begin{picture}(20,20)(0,8)
\thicklines
\put(0, 0){\line(1,0){10}}
\put(0,10){\line(1,0){20}}
\put(0,20){\line(1,0){20}}
\put( 0, 0){\line(0,1){20}}
\put(10, 0){\line(0,1){20}}
\put(20,10){\line(0,1){10}}
\put( 2,11.8){\bf 1}
\put(12,11.8){\bf 2}
\put( 2, 1.8){\bf 3}
\end{picture}
\ =\
\begin{picture}(20,10)(0,3)
\thicklines
\put(0, 0){\line(1,0){20}}
\put(0,10){\line(1,0){20}}
\put( 0, 0){\line(0,1){10}}
\put(10, 0){\line(0,1){10}}
\put(20, 0){\line(0,1){10}}
\put( 2, 1.8){\bf 1}
\put(12, 1.8){\bf 2}
\end{picture}
\circ
\setlength{\unitlength}{1.3pt}%
\begin{picture}(10,10)(0,3)
\thicklines
\put(0, 0){\line(1,0){10}}
\put(0,10){\line(1,0){10}}
\put( 0, 0){\line(0,1){10}}
\put(10, 0){\line(0,1){10}}
\put( 2, 1.8){\bf 1}
\end{picture}
\ ,\raisebox{-20pt}{\rule{0pt}{5pt}}
$$
so $\circ$ is non-commutative.
To show associativity,
let $R$, $S$, and $T$ be tableaux.
Then
$$
R\circ (S\circ T) = \sum R \bigcup W,
$$
the sum over $W$ Knuth equivalent to $S\bigcup V$,
where $V$ is Knuth equivalent to $T$.
Let $U'$ be the first $|S|$ entries in $W$, and
$V'$ the last $|T|$ entries, each decreased by $|S|$,
thus,
$$
R\circ (S\circ T) = \sum R \bigcup U' \bigcup V',
$$
the sum over  $U'$  Knuth equivalent
to $S$ and $V'$ to $T$, which is $(R\circ S)\circ T$.
\QED

Let $m<n$.
For a tableau $T$ of shape $\lambda$, let
$\phi(T)$ be the Schur function $s_\lambda$.
Define $\phi_{m,\,n}(T)$ to be 0 if $\lambda_1+m \geq n$
or $\lambda_{m+1}\neq 0$ and $\sigma_\lambda$ otherwise.
Then $\phi$ and $\phi_{m,n}$ are, respectively, additive surjections
from $\Lambda$
to the algebra of symmetric functions and to $A^*\G_m{\bf P}^n$.
\medskip

\noindent{\bf Theorem G.} {\em
The maps $\phi$ and $\phi_{m,\,n}$ are $\bf Z$-algebra
homomorphisms.
}
\medskip

\noindent{\bf Proof:}
For any tableaux $S$ and $T$ of shape $\nu$ and partitions $\lambda$
and $\mu$, there is a natural bijection
(given by dual equivalence of Haiman~\cite{Haiman_dual_equivalence})
between the set of
tableaux with shape $\lambda/\mu$ Knuth equivalent to $S$ and those
Knuth equivalent to $T$, and this common number is
$c^\lambda_{\mu\,\nu}$.
This shows that $\phi$ is an algebra homomorphism.
It follows that
$\phi_{m,\,n}$ is as well.
\QED

\section{Enumerative Geometry and Arrangements Over Finite Fields}

A main result of this paper, Theorem~C, shows that any Schubert-type
enumerative problem concerning lines in projective space may be
solved over $\bf R$.
By `solved' over a field $k$, we mean
there are flags in ${\bf P}^n_k$ determining Schubert varieties which meet
transversally in finitely many points, all of which are
defined over $k$.
Given an enumerative problem, we feel it is legitimate
to inquire over which (finite) fields it may be solved.

We present two families of enumerative problems
for which this question may be resolved, and consider the problem
of finding arrangements over finite fields.

\subsection{The $n$ lines meeting four $(n-1)$-planes in ${\bf P}^{2n-1}$.}
Given three  non-intersecting $(n-1)$-planes $L_1,L_2,$ and $L_3$
in ${\bf P}^{2n-1}$,
there are coordinates $x_1,\ldots,x_{2n}$
such that
\begin{eqnarray*}
L_1 &:& x_1 = x_2 = \cdots = x_n = 0\\
L_2 &:& x_{n+1} = \cdots = x_{2n} = 0\\
L_3 &:& x_1 - x_{n+1} = \cdots = x_n - x_{2n}= 0
\end{eqnarray*}
One may  check that
$ \Omega_{L_1}\bigcap\Omega_{L_2}\bigcap\Omega_{L_3}$
is a transverse intersection, and if $\Sigma_{1,n-1}\subset {\bf P}^{2n-1}$
is the union of the lines meeting each of
 $L_1,L_2,$ and $L_3$, then $\Sigma_{1,n-1}$ is the image of the standard
Segre embedding of
${\bf P}^1\times{\bf P}^{n-1}$ into ${\bf P}^{2n-1}$
(cf.~\cite{Harris_geometry}):
$$
\psi:[a,b]\times[y_1,\ldots,y_n] \longmapsto
[ay_1,\ldots,ay_n,by_1,\ldots,by_n].
$$
The lines meeting $L_1$, $L_2$, and $L_3$ are
the images of ${\bf P}^1\times \{p\}$, for $p\in{\bf P}^{n-1}$.
$\Sigma_{1,n-1}$ has degree $n$, so a general $(n-1)$-plane
$L_4$ meets  $\Sigma_{1,n-1}$ in $n$ distinct points, each determining
a line meeting $L_1,\ldots,L_4$.
These lines, $\ell_1,\ldots,\ell_n$, meet $L_1$ in distinct points
which span $L_1$.
Changing coordinates if necessary, we may
assume $\ell_j$ is the span of $x_j$ and $x_{n+j}$.

For $1\leq j\leq n$, let $p_j = [\alpha_j,\beta_j]\in {\bf P}^1$ be the
first coordinate of $\psi^{-1}(\ell_j\cap L_4)$.
Then
$$
L_4 \quad :\quad \beta_1x_1 -\alpha_1 x_{n+1}=\cdots=
 \beta_nx_n - \alpha_n x_{2n}=0.
$$
Also, $p_1,\ldots,p_n$ are distinct; otherwise
$L_4\cap \Sigma_{1,n-1}$ contains a line.

Thus, if $k$ has at least $n-1$ elements, this enumerative problem
may be solved over $k$.

\subsection{The $n$ lines meeting a fixed line and $n+1$
$(n-1)$-planes in ${\bf P}^{n+1}$.}
A line $\ell$ and $(n-1)$-planes
$K_1,\ldots,K_n$ in ${\bf P}^{n+1}$
are independent if for every $p\in \ell$, the hyperplanes
$\Gamma_i(p) = \Span{p,K_i}$, for $1\leq i\leq n$, meet in a line.
In this case, the union
$$
S_{1,n-1} = \bigcup_{p\in \ell}
\Gamma_1(p)\cap \cdots\cap \Gamma_n(p)
$$
is a rational normal surface scroll.
Moreover, the lines meeting each of $\ell$, $K_1,\ldots,K_n$
are precisely those lines
$\lambda(p) = \Gamma_1(p)\cap \cdots\cap \Gamma_n(p)$
for $p\in \ell$.

Since $S_{1,n-1}$ has degree $n$, a general $(n-1)$-plane
$K_{n+1}$ meets $S_{1,n-1}$ in $n$ distinct points, each determining a
line $\lambda(p)$ which meets
$\ell,K_1,\ldots,K_{i+1}$.
If $k$ is finite with $q$ elements,
there are only $q+1$ lines $\lambda(p)$ defined over $k$.
Thus it is necessary that $q\geq n-1$ to solve this problem over $k$.
We show this condition suffices.

All rational normal surface scrolls are projectively equivalent,
(cf.~\cite{Harris_geometry}, \S9),
thus we may assume that $S_{1,n-1}$ has the following standard form.
Let $x_1,x_2,y_1,\ldots,y_n$ be coordinates for ${\bf P}^{n+1}$
where $\ell$ has equation $y_1 = \cdots = y_n = 0$.
Then for $p = [a,b,0,\ldots,0] \in \ell$,
$\lambda(p)$ is the linear span of $p$ and the point
$[0,0,a^{n-1},a^{n-2}b,\ldots,ab^{n-2},b^{n-1}]$.

Let $\alpha_1,\ldots,\alpha_n\in {\bf P}^1$ be distinct points.
Let $F = \sum_{i=0}^n A_i b^i a^{n-i}$
be a form on ${\bf P}^1$ vanishing at $\alpha_1,\ldots,\alpha_n$.
Define $K_{n+1}$ by the vanishing of the two linear forms
$$
\Lambda_1 \ :\  x_2 - y_1\ \ \ \ \ \
\Lambda_2 \ :\  A_0 x_1 + A_1y_1 + \cdots + A_n y_n.
$$
The intersection of $S_{1,n-1}$ and the hyperplane defined by
$\Lambda_1$ is the rational normal curve
$$
\psi : [a,b] \longmapsto [a^n,a^{n-1}b,a^{n-1}b,
a^{n-2}b^2,\ldots,a b^{n-1},b^n].
$$
Since $\psi^* (\Lambda_2) = F$, the lines meeting each of
$\ell,K_1,\ldots,K_{n+1}$ are
$\lambda(\alpha_1),\ldots,\lambda(\alpha_n)$.

Thus, if $k$ has at least $n-1$ elements, this enumerative problem
may be solved over $k$.
\medskip

These  two families
are the only non-trivial examples of Schubert-type enumerative
problems for which we know an explicit description
of their solutions.
Each of these problems can be solved over any field $k$
where $\#{\bf P}^1_k$ exceeds the number of solutions.
It would be interesting to find explicit
solutions to other enumerative problems to test whether this holds more
generally.

\subsection{Arrangements over Finite Fields}
In \S\ref{sec:arrangements} we remarked it is possible to construct
arrangements over some
finite fields.
Here we  show how.

Recall that an arrangement is complete flag $\Fdot$
and $2n-3$ hyperplanes $H_2,\ldots,H_{2n-2}$
such that for any $p$,
\begin{enumerate}
\item  $H_p \cap  F_{\hfl{p}} =  F_{\hfl{p}+1}$,
\item  For tableaux $S, T$ of degrees at most $p-1$, if
$H_T\cap H_S \subset H_p$, then $H_T\cap H_S = F_{\hfl{p}+1}$.
\end{enumerate}
We give an equivalent set of conditions.
For any subset $A\neq \emptyset$ of $\{2,3,\ldots,2n-2\}$, let $H_A$ be
$\bigcap\{H_i\,|\, i\in A\}$.
Set $H_{\emptyset} = \Pn$.
\medskip

\noindent{\bf Lemma.}\label{lemma:alt_def_arrangement}
{\em A complete flag $\Fdot$ and hyperplanes $\,H_2,\ldots,H_{2n}$
constitute an arrangement if and only if
for each $m = 1,\ldots,n-1$, they satisfy
\begin{enumerate}
\item [1$'$.]  $H_{2m}\cap F_m = F_{m+1}$ and $F_{m+1}\subset H_{2m+1}$.
\item [2$'$.]  For any $A\subset \{2,3,\ldots,2m\}$ where
$H_A$ has codimension $\#A \leq m$, $H_A\not\subset H_{2m+1}$.
\end{enumerate}
}
\medskip

\noindent{\bf Proof:}
Let $H_2,\ldots,H_{2n-2}$ and $\Fdot$  satisfy 1$'$ and 2$'$.
We show they constitute an arrangement by induction on $m$.
Suppose that for $p<2m$ conditions 1 and 2 for arrangements are satisfied.
We show that 1 and 2 are satisfied for $p = 2m$ and $2m+1$.
For $p=2m$, 1 and 1$'$ are equivalent.
Moreover, if $S,T$ are tableaux of degree less than $2m$,
then $F_m \subset H_T\cap H_S$, so $H_T\cap H_S\not\subset H_{2m}$,
so 2 is satisfied.

Thus $H_2,\ldots,H_{2m}$, $F_1,\ldots,F_{m+1}$
constitute an arrangement in ${\bf P}^n/F_{m+2}\simeq {\bf P}^{m+1}$.
Then by Lemma~\ref{lemma:arrangements}, if $T$ has shape $(\alpha,\beta)$
with $\alpha+\beta < 2m$, $H_T$ has codimension $\beta$ and is not equal to
$F_\beta$.

Let $S,T$ be  tableaux of degree at most $2m$
and suppose $F_{m+1}\neq H_S\cap H_T$.
Let $B$ be the union of the second rows of $S$ and $T$.
Since $H_B = H_S\cap H_T$ has codimension $s \leq m$, there is a
set $A\subset B$ of order $s$ with $H_A = H_B$.
By 2$'$, $ H_S\cap H_T = H_A \not\subset H_{2m+1}$.
So $H_{2m+1}$ satisfies 1 and 2.

Conversely, suppose $H_2,\ldots,H_{2n-2}$ and  $\Fdot$
constitute an arrangement.
These satisfy 1$'$.
To show they satisfy 2$'$,
let $A\subset \{2,\ldots,2m\}$ where $H_A$ has
codimension $s = \# A$.
Suppose $A = a_1<\cdots<a_s$ and let $j$ be the largest index such that
$a_j<2j$.
If $j=0$, then $A$ is the second row of a tableau $T$,
so $H_A = H_T \not\subset H_{2m+1}$.
If $j\neq 0$, then $j\geq 2$.
Since $p<2j$ implies $F_j\subset H_p$ and $H_A$ has codimension $s$, we
have $F_j = H_{a_1}\cap \cdots\cap H_{a_j}$.
An induction using condition 1$'$ shows that
$F_j = H_3\cap H_2\cap H_4\cap \cdots\cap H_{2j-2}$.
Thus $H_A = H_{A'}$ where
$$
A' = 3,2<\cdots<2j-2<a_{j+1}<\cdots<a_s,
$$
and if $i>j$, $a_i\geq 2i$.
Let $B: 2<\cdots<2j-2<a_{j+1}<\cdots<a_s$.
then we see that $B$ is the second row a tableau $T$ of degree at most $2s$.
Let $S$ be the tableau of degree 3 whose second row consists of $3$.
Since $s\leq m$, $F_{m+1} \neq H_S\cap H_T$,
so $H_A = H_S\cap H_T \not\subset H_{2m+1}$.
\QED

\subsection{Corollary.}\label{sec:refined_condition} {\em
Condition $2'$ may be replaced by
\begin{enumerate}
\item[2$''$.] If $A: 3,2<\cdots<2j-2<a_{j+1}<\cdots<a_m\leq 2m$
satisfies
$i>j$ implies $a_i\geq 2i$, then $H_A\not\subset H_{2m+1}$.
\end{enumerate}}
\medskip

\noindent{\bf Proof:}
Suppose $H_2,\ldots,H_{2n-2}$ and $\Fdot$ satisfy 1$'$ and 2$'$
and $A : 3,2<\cdots<2j-2<a_{j+1}<\cdots<a_m\leq 2m$ satisfies
$i>j$ implies $a_i \geq 2i$.
Then $2<\cdots<2j-2<a_{j+1}<\cdots<a_m$ is the second row of
a tableau $T$.
By Lemma~\ref{lemma:arrangements}, $H_{T}$ has codimension $m-1$ and
does not equal $F_{m-1}$.
Thus $H_T\not\subset H_3$ so $H_A = H_T\cap H_3$ has codimension $m$.

Conversely, suppose
$H_2,\ldots,H_{2n-2}$ and $\Fdot$ satisfy 1$'$ and 2$''$.
{}From the proof of the previous Lemma, it suffices to know
2$'$ for those subsets $A$ of the form
$3,2<\cdots<2j-2<a_{j+1}<\cdots<2s\leq 2m$, where $i>j$ implies $a_i\geq 2i$.
If $B = A\cup\{2s+2,\ldots,2m\}$, then $B$ also has the form
in 2$''$.
So  $H_B \not\subset H_{2m+1}$.
Since $H_B\subset H_A$, $H_A\not\subset H_{2m+1}$,
showing 2$'$ holds for $H_2,\ldots,H_{2n-2},\Fdot$.
\QED

We  estimate the size of a field $k$ necessary to
construct an arrangement.

\subsection{Theorem.}\label{thm:arrangement_finite_field}
{\em
There exists an arrangement in  ${\bf P}^n_k$ if the order of $k$
is at least
$$
\frac{(2n-4)!}{(n-2)!(n-1)!}\ +\ \sum_{i=1}^{n-4} \frac{(2i)!}{i! (i+1)!}.
$$
}
\medskip

\noindent{\bf Proof:}
Consider the problem of inductively constructing an arrangement
in ${\bf P}^n_k$ satisfying 1$'$ and 2$''$ of
\S\S\ref{lemma:alt_def_arrangement} and~\ref{sec:refined_condition}.
Since it is always possible to find a hyperplane not
containing any particular proper linear subspace of $\Pn$, the only possible
obstruction is the selection of hyperplanes $H_{2m+1}\supset F_{m+1}$
satisfying 2$''$ for $m=0,1,\ldots,n-2$.

Let $\check{\bf P}^m$ be the set of hyperplanes defined over $k$
containing $F_{m+1}$.
Every codimension $m$ subspace $H_A$ containing $F_{m+1}$
determines a hyperplane $\check{H}_A$ in $\check{\bf P}^m$
consisting of those hyperplanes $H$ of ${\bf P}^n$ containing
$H_A$.
Thus there exists a hyperplane $H_{2m+1}$ satisfying 2$'$ if,
as sets of $k$-points,
$$
X = \check{\bf P}^m - \bigcup_{A\in {\cal S}} \check{H}_A\neq \emptyset,
$$
Where ${\cal S}$ is the set of those sequences
$3,2<\cdots<2j-2<a_{j+1}<\cdots<a_m\leq 2m$ such that
if $i>j$, then $a_i\geq 2i$.
We claim
$\#{\cal S} = s_m =  \frac{(2m)!}{m!(m+1)!} +
\sum_{i=1}^{m-2} \frac{(2i)!}{i! (i+1)!}$.
Suppose $k$ has $q\geq s_{n-2}$ elements.
Since $\check{\bf P}^m$ has $(q^{m+1}-1)/(q-1)$ elements and
each $\check{H}_A$  has $(q^m-1)/(q-1)$ elements, $X$
is non-empty if
$q^{m+1}-1 > (q^m-1)s_m$.
This holds as
$$
\left\lfloor\frac{q^{m+1}-1}{q^m-1}\right\rfloor \geq q \geq s_{n-2}\geq s_m.
$$

To enumerate ${\cal S}$, let
$\{3,2<\cdots<2j-2<a_{j+1}<\cdots<a_m\} \in {\cal S}$.
If $b_j = a_{j+i} - 2j$, then
$b_1,\ldots, b_{m-j}$ is the second row of a tableau of shape
$(m-j,\,m-j)$.
Conversely, if $b_1,\ldots, b_{m-j}$ is the second row of a
tableau of shape $(m-j,\,m-j)$, then
$$
\{3, 2<\cdots<2j-2<b_1+2j<\cdots<b_{m-j}+2j\} \in {\cal S}.
$$
Let ${\cal T}_s$ be the set of tableaux of shape $(s,s)$.
These considerations show there is a bijection
$$
{\cal S} \longleftrightarrow
{\cal T}_m\cup{\cal T}_{m-2}\cup{\cal T}_{m-3}\cup\ldots\cup{\cal T}_0.
$$
Noting that $\#{\cal T}_s = \frac{(2s)!}{s!(s+1)!}$,
by the hook length formula of Frame, Robinson, and Thrall~\cite{FRT},
shows that the order of ${\cal S}$ is
$\frac{(2m)!}{m!(m+1)!}+\sum_{i=1}^{m-2} \frac{(2i)!}{i! (i+1)!}$.
\QED

This result is not the best possible:
For ${\bf P}^4$, this gives $q\geq 5$, but
arrangements in ${\bf P}^4$ may be constructed over
the field with three elements.

\end{document}